\documentclass[graybox]{svmult}

\usepackage{mathptmx}       
\usepackage{helvet}         
\usepackage{courier}        
\usepackage{type1cm}        
\usepackage{dsfont}
\usepackage{makeidx}         
\usepackage{graphicx}        
\usepackage{multicol}        
\usepackage[bottom]{footmisc}

\makeindex

\begin{document}

\title*{Two-dimensional fluid queues with temporary assistance}
\author{Guy Latouche, Giang T. Nguyen, Zbigniew Palmowski}
\institute{Guy Latouche \at Universit\'{e} Libre de Bruxelles, D\'{e}partement d'Informatique, CP 212, Boulevard du Triomphe, 1050 Brussels, Belgium, \email{latouche@ulb.ac.be}
\and Giang T. Nguyen  \at Universit\'{e} Libre de Bruxelles, D\'{e}partement d'Informatique, CP 212, Boulevard du Triomphe, 1050 Brussels, Belgium, \email{giang.nguyen@ulb.ac.be}
\and Zbigniew Palmowski \at University of Wroc{\l}aw, Mathematical Institute, 50-384 Wroc{\l}aw, Poland, \email{zbigniew.palmowski@gmail.com}}

\maketitle

\abstract{We consider a two-dimensional stochastic fluid model with $N$ ON-OFF inputs and temporary assistance, which is an extension of the same model with $N = 1$ in Mahabhashyam \emph{et al.} (2008). The rates of change of both buffers are piecewise constant and dependent on the underlying Markovian phase of the model, and the rates of change for Buffer~2 are also dependent on the specific level of Buffer~1. This is because both buffers share a fixed output capacity, the precise proportion of which depends on Buffer~1. The generalization of the number of ON-OFF inputs necessitates modifications in the original rules of output-capacity sharing from Mahabhashyam \emph{et al.} (2008) and considerably complicates both the theoretical analysis and the numerical computation of various performance measures. \\
\\
We derive the marginal probability distribution of Buffer~1, and bounds for that of Buffer~2. Furthermore, restricting Buffer~1 to a finite size, we determine its marginal probability distribution in the specific case of $N = 1$, thus providing numerical comparisons to the corresponding results in Mahabhashyam \emph{et al.} (2008) where Buffer~1 is assumed to be infinite. }

\section{Introduction}
Stochastic fluid models have a wide range of applications, such as water reservoir operational control, industrial and computer engineering, risk analysis, environmental analysis and
telecommunications. In particular, they have been used in telecommunication modeling since the seminal paper~\cite{anick82}. With the advent of differentiated services, buffers have in a very natural way
become multidimentional. To give another example, that of decentralized mobile networks, callers transmit data via each other's equipment and it is necessary to determine the appropriate fractions
of a caller capacity, be it buffer space or power, that may be allocated to other users. \\
\\
In computer processing, a situation where the problem of effective resource-sharing can arise is when there are more tasks than schedulers that can process them. Aggarwal \emph{et al.}~\cite{aggarwal59} consider this problem in the particular setting of two ON-OFF streams of tasks: routine and non-routine, and one Central Processing Unit (CPU) to serve both streams, which one, specifically, being determined by a workload threshold. The CPU serves routine or non-routine tasks, depending on whether the amount of workload for the routine tasks is above or below the threshold, respectively. In order to determine the optimal threshold value that minimizes the weighted sum of the probability of exceeding undesirable workload limits, the authors derive the workload distribution of routine tasks, and approximate that of non-routine tasks. Mahabhashyam \emph{et al.}~\cite{resourcesharing} extends the resource-sharing model to allow a partial split of the CPU's capacity. More specifically, the CPU serves routine tasks when their accumulated workload is above the threshold; the CPU serves, according to some predetermined proportion, both routine and non-routine tasks when the threshold is not exceeded; and the CPU serves non-routine tasks when there is no routine task left. \\
\\
We generalize this model further, by allowing the input model to
better fit an environment where multiple users independently
decide when to use the system, thereby allowing for the intensity of
the load to vary in time.  Specifically, each input stream of fluid
is formed by $N$ exponential ON-OFF sources, with $N \geq 1$, and we
analyze the model using a two-dimensional stochastic fluid. \\
\\
A Markov-modulated single-buffer fluid model is a two-dimensional Markov process $\{X(t),\varphi(t): t \in \mathds{R}^{+}\}$, where $X(t)$ is the continuous \emph{level} of the buffer, and $\varphi(t)$ is the discrete \emph{phase} of the underlying irreducible Markov chain that governs the rates of change. A practical and well-studied case is \emph{piecewise constant} rates: the fluid is assumed to have a constant rate $c_i$ when $\varphi(t) = i$, for $i$ in a finite state space $\mathcal{S}$. The traditional approach for obtaining performance measures of Markov-modulated single-buffer fluids with piecewise constant rates is to use spectral analysis (see, among others, \cite{kosten74, anick82, mitra88, stern91, elwalid93}). Over the last two decades, matrix analytic methods have gained a lot of attention as an alternative and algorithmically effective approach for analyzing these standard fluids (see, for instance, \cite{ram99, akar04, unboundedfluids, anafinite, bo08, algosforlst,leveldep, boundedmodel}).\\
\\
The mathematical model we consider is a Markov process $\{X(t), Y(t),$ $ \varphi_1(t), \varphi_2(t): t \in \mathds{R}^{+}\}$, where $X(t) \geq 0$ and $Y(t) \geq 0$ represent the levels of Buffers 1 and 2, respectively. At a given time $t \geq 0$, the rates of change of Buffer~1 depend only on the underlying Markovian phase $\varphi_1(t)$; the rates of change of Buffer~2, on the other hand, depend on both $\varphi_2(t)$ and $X(t)$. This is because while each buffer receives its own input sources, both share a fixed output capacity $c$, in proportion dependent on the level of Buffer~1. More specifically, Buffer~$j$ receives $N$ ON-OFF input sources, each has exponentially distributed ON- and OFF- intervals at corresponding rates $\alpha_j$ and $\beta_j$, and continuously generates fluid at rate $R_j$ during ON- intervals, for $j = 1, 2$. When the fluid level $X(t)$ of Buffer~1 is above the threshold $x^* > 0$, Buffer~1 is allocated the total shared output capacity $c$, leaving Buffer~2 without any; when $0 < X(t) < x^*$, Buffer~$j$ has output capacity $c_j$, $c_1 + c_2 = c$; and when $X(t) = 0$, Buffer~1 has output capacity $\min\{i R_1, c_1\}$, and Buffer~2 $\;c - \min\{i R_1,c_1\}$, where $i$ is the number of inputs of Buffer~1 being on at the time $t$. \\
\\
The generalization of the number of ON-OFF inputs necessitates modifications in the original rules of output-capacity sharing from Mahabhashyam \emph{et al.} \cite{resourcesharing}. When $X(t) = 0$, the policy in the single ON-OFF input model is to allocate the total capacity $c$ to Buffer~2. The totality rule is logical when there is only one ON-OFF input for each buffer: Buffer~1 is empty only when its input is off, and in that case, Buffer~2 can receive the whole output capacity $c$, until the moment the input of Buffer~1 is on again. Here, it is possible for Buffer~1 to be empty while $i$ inputs are on, for $0 < i \leq \lfloor \frac{c_1}{R_1} \rfloor$. Under these circumstances, assigning the total output capacity $c$ to Buffer~2 would immediately cause Buffer~1 to try to increase from level 0, consequently grabbing back $c_1$ amount of output capacity. However, as $i \leq \lfloor \frac{c_1}{R_1} \rfloor$, the output capacity $c_1$ would be sufficient to empty Buffer~1, forcing it to give away the whole output capacity $c$ to Buffer~2, etc. Therefore, applying the original totality rule at $X(t) = 0$ for the generalized $N$ ON-OFF input model would potentially lead to inconsistency. \\
\\
The behavior described above at level 0 for Buffer~1 when $0 \leq i \leq \lfloor \frac{c_1}{R_1} \rfloor$ is referred to as being \emph{sticky} \cite{leveldep}, a property arising when net rates of the buffer for the same Markovian phase but different levels are different in a particular way that makes it unable to go up or down, thus remaining stuck at a level until the background Markov chain switches to a non-sticky phase. In our model, by allocating $i R_1$ output capacity to Buffer~1 and $c - i R_1$ to Buffer~2 when $X(t) = 0$ and $0 \leq i \leq \lfloor \frac{c_1}{R_1} \rfloor$, we let Buffer~1 remain at level zero, while eliminating potential uncertainty and utilizing the total output capacity in the most effective way. For the same reason, when $X(t) = x^*$ and $\lceil \frac{c_1}{R_1} \rceil \leq i \leq \lfloor \frac{c}{R_1} \rfloor$, the output capacity is $i R_1$ for Buffer~1, and $c - i R_1$ for Buffer~2. While the stickiness, borne in the generalization of the number of ON-OFF inputs, necessitates only slight modifications in the output-capacity allocation policy, it considerably complicates the analysis and numerical computation of performance measures of the model. To deal with this complication, we employ a mixture of tools from both dominant approaches: spectral analysis and matrix analytic methods. \\
\\
One may change the system in many ways and still use the same method. For example, in the last part of this paper, we restrict Buffer~1 to a finite size, but keep Buffer~2 being infinite. This affects the analysis of Buffer~1, but the analytical expressions for Buffer~2 remain unchanged. We take $N = 1$ there for better illustration.  \\
\\
The rest of the paper is organized as follows: in Section~\ref{sec:refmodel}, we formulate the model mathematically. Assuming that both buffer sizes are infinite, we derive the marginal probability distribution of Buffer~1 in Section~\ref{sec:infb1}, and bounds for those of Buffer~2 in Section~\ref{sec:infb2}. In Section~\ref{sec:finb1},
restricting Buffer~1 to a finite size, we determine its marginal probability distribution in the particular case of $N = 1$, thus providing numerical comparisons to the corresponding results in \cite{resourcesharing}, where Buffer~1 is assumed to be infinite.

\section{Reference model} \label{sec:refmodel}
Consider a four-dimensional Markov process $\{X(t), Y(t),$ $\varphi_1(t), \varphi_2(t): t \in \mathds{R}^{+}\}$, where $X(t) \geq 0$ and $Y(t) \geq 0$ are the levels in Buffers~1 and 2, respectively, and for $j = 1,2$, $\varphi_j(t)$ represents the phase of the background irreducible Markov chain for Buffer~$j$ with finite state space $\mathcal{S} = \{0, \ldots, N\}$ with $N \geq 1$; state $i \in \mathcal{S}$ indicates that $i$ ON-OFF inputs are on. The generator $T_j$ for $\{\varphi_j(t)\}$ is
\begin{eqnarray*}
T_j & = & \left[\begin{array}{ccccc}
           * &    N\beta_j &                   &                         & \\
    \alpha_j &           * & (N - 1)\beta_j    &                         & \\
             & \ddots      &   \ddots          & \ddots                  & \\
             &             &   (N - 1)\alpha_j &                       * & \beta_j \\
             &             &                   & N\alpha_j               & *
\end{array}\right],
\end{eqnarray*}
with each diagonal element $*$ defined appropriately such that each row sum of $T_j$ is 0. For $i_1, i_2 \in \mathcal{S}$, we denote by $\dot{x}_{i_1}$ and $\dot{y}_{i_2}$ the respective net rates for Buffer~1 in phase $i_1$ and Buffer~2 in phase $i_2$. For $X(t) > x^*$ and $Y(t) > 0$,
\begin{eqnarray*}
         \dot{x}_{i_1} & = & i_1 R_1 - c, \\
         \dot{y}_{i_2} & = & i_2 R_2;
\end{eqnarray*}
for $X(t) = x^*$ and $Y(t) > 0$,
\begin{eqnarray*}
         \dot{x}_{i_1} & = & 0 \quad \mbox{ for } \lceil \frac{c_1}{R_1} \rceil \leq i_1 \leq \lfloor \frac{c}{R_1} \rfloor,\\
                       & = & i_1 R_1 - c_1 \quad \mbox{ otherwise}, \\
         \dot{y}_{i_2} & = & i_2 R_2 - (c - i_1R_1) \quad \mbox{ for } \lceil \frac{c_1}{R_1} \rceil \leq i_1 \leq \lfloor \frac{c}{R_1} \rfloor,\\
                       & = & i_2 R_2 - c_2 \quad \mbox{ otherwise};
\end{eqnarray*}
for $0 < X(t) < x^*$ and $Y(t) > 0$,
\begin{eqnarray*}
         \dot{x}_{i_1} & = & i_1 R_1 - c_1,  \\
         \dot{y}_{i_2} & = & i_2 R_2 - c_2;
\end{eqnarray*}
and for $X(t) = 0$ and $Y(t) > 0$:
\begin{eqnarray*}
\dot{x}_{i_1} & = & 0 \quad \mbox{ for } 0 \leq i_1 \leq \lfloor \frac{c_1}{R_1} \rfloor,\\
                    & = & i_1 R_1 - c_1 \quad \mbox{ otherwise}, \\
\dot{y}_{i_2} & = & i_2 R_2 - (c - i_1 R_1) \quad \mbox{ for } 0 \leq i_1 \leq \lfloor \frac{c_1}{R_1} \rfloor,\\
                    & = & i_2 R_2 - c_2 \quad \mbox{ otherwise}.
\end{eqnarray*}
For $Y(t) = 0$, $\dot{y}_{i_2}$ is the maximum between 0 and the net rate of Buffer~2 in $i_2 \in \mathcal{S}$ when $Y(t) > 0$. \\
\\
We assume that $NR_j > c$, $\frac{c_j}{R_j}, \frac{c}{R_j} \not\in \mathds{N}$, and the system is positive recurrent. The first assumption ensures that for $X(t) > x^*$, the set of states for which the net rates of Buffer~1 are positive is non-empty. We impose the second assumption to avoid having states with zero rates for Buffer~1 when $X(t) \notin \{0, x^*\}$, and for Buffer~2 when $Y(t) \neq 0$. This assumption is purely to simplify some technical details, without any loss of generality, as any single-buffer fluid model with zero rates can be transformed into a single-buffer fluid model without zero rates \cite{ana02}. The third assumption is equivalent to
\begin{eqnarray} \label{eqn:cond}
\sum_{i = 0}^{N} i R_1 q^{(1)}_i  + \sum_{i = 0}^{N} i R_2 q^{(2)}_i < c_1 + c_2,
\end{eqnarray}
where $\underline{q}^{(1)}$ and $\underline{q}^{(2)}$ are the stationary probability vectors of $T_1$ and $T_2$, respectively. The inequality $(\ref{eqn:cond})$ is obvious when considering the stability condition for the equivalent single-buffer fluid model with a constant output $c_1 + c_2$ and $2N$ exponential ON-OFF inputs, half of which switch on at rate $\beta_1$ and switch off at rate $\alpha_1$, and the other half switch on at rate $\beta_2$ and switch off at rate $\alpha_2$. For $i = 1, 2$ and for $n = 1, \ldots, N$,
\begin{eqnarray*}
q_0^{(i)} & = &\frac{\alpha_i^N}{(\alpha_i + \beta_i)^N}, \\
q_n^{(i)} & = & q_0^{(i)} {{N}\choose{n}}(\beta_i/\alpha_i)^{n},
\end{eqnarray*}
which reduces (\ref{eqn:cond}) to
\begin{eqnarray}
N\frac{R_1\beta_1}{\alpha_1 + \beta_1} + N\frac{R_2\beta_2}{\alpha_2 + \beta_2} < c_1 + c_2.
\end{eqnarray}
\section{Infinite Buffer~1} \label{sec:infb1}
To analyze Buffer~1 when $N = 1$, Mahabhashyam \emph{et al.} \cite{resourcesharing} consider an equivalent system of two standard single sub-buffers, each with a single ON-OFF input, one sub-buffer with constant output capacity $c_1$ and the other with constant output capacity $c$. Decomposing Buffer~1 in this fashion, the authors show that the marginal probability distribution of Buffer~1 can be obtained by appropriately combining the average time of going up from $x^*$ and then going down to $x^*$ in Sub-buffer 1, and the average time of going down from $x^*$ and then going up to $x^*$ in Sub-buffer 2. The authors determine analytic expressions for the former average time by using, from \cite{narayanan96}, the busy period distribution of a standard single buffer with one exponential ON-OFF input and constant output capacity, and for the latter by establishing a pair of partial differential equations, transferred into ordinary differential equations and then solved by a spectral decomposition technique. \\
\\
In this paper, for general $N \geq 1$, we analyze Buffer~1 by applying matrix analytic methods. With this approach, while it is not simple to obtain closed-form expressions for $N \geq 2$, we can obtain various performance measures numerically using fast convergent algorithms (see, most relevantly, \cite{bo08,leveldep} and the references therein). The focus of this section is the marginal probability distribution for Buffer~1. \\
\\
We refer to $X(t) = 0$ and $X(t) = x^*$ as \emph{boundaries} $\circ$ and $*$, and $0 < X(t) < x^*$ and $X(t) > x^*$ as \emph{bands} 1 and 2. While $T_1$ governs the transitions of $\{\varphi_1(t)\}$ for all $X(t) \geq 0$, the rate of Buffer~1 in the same phase varies between boundaries and bands. Therefore, we partition $\mathcal{S}$ differently for each boundary and each band. We denote, respectively, by $\mathcal{S}^{(\bullet)}_d$, $\mathcal{S}^{(\bullet)}_s$ and $\mathcal{S}^{(\bullet)}_u$ the sets of states with negative, zero and positive net rates when Buffer~1 is at boundary $\bullet$, for $\bullet \in \{\circ, *\}$, and by $\mathcal{S}^{(k)}_{-}$ and $\mathcal{S}^{(k)}_{+}$ the sets of states with negative and positive net rates when Buffer~1 is in band $k$, for $k = 1, 2$. Then,
$ \mathcal{S}  = \mathcal{S}_{s}^{(\circ)} \cup \mathcal{S}_{u}^{(\circ)} = \mathcal{S}^{(1)}_{-} \cup \mathcal{S}^{(1)}_{+}  = \mathcal{S}_{d}^{(*)} \cup \mathcal{S}_{s}^{(*)} \cup \mathcal{S}_{u}^{(*)} = \mathcal{S}^{(2)}_{-} \cup \mathcal{S}^{(2)}_{+}, $
with
\begin{eqnarray*}
\mathcal{S}^{(\circ)}_s & = & \mathcal{S}^{(1)}_{-} =  \mathcal{S}^{(*)}_d = \mathcal{S}_{-}^{(2)} = \{0,\ldots, \lfloor \frac{c_1}{R_1} \rfloor\}, \\
\mathcal{S}^{(\circ)}_u & = & \mathcal{S}^{(1)}_{+} = \mathcal{S}_{+}^{(2)} = \{\lceil \frac{c_1}{R_1} \rceil, \ldots, N\}, \\
\mathcal{S}_s^{(*)} & = & \{\lceil \frac{c_1}{R_1}\rceil, \ldots, \lfloor \frac{c}{R_1}\rfloor\}, \mathcal{S}^{(*)}_u = \{\lceil \frac{c}{R_1} \rceil,\ldots, N\}.
\end{eqnarray*}
For each band $k$, we partition $T_1$ into sub-matrices $T^{(k)}_{\ell m}$, of which each element $[T^{(k)}_{\ell m}]_{ij}$ is the transition rate from $i \in \mathcal{S}^{(k)}_{\ell}$ to $j \in \mathcal{S}^{(k)}_{m}$, and we denote by $C_{\ell}^{(k)}$ the diagonal matrix of absolute net rates for $i \in \mathcal{S}_{\ell}^{(k)}$:
\begin{eqnarray*}
C_{-}^{(1)} & = &  \left[\begin{array}{cccc}
|- c_1| &            &        & \\
          & |R_1 - c_1| &        & \\
          &              & \ddots & \\
          &              &        & |\lfloor \frac{c_1}{R_1} \rfloor R_1 - c_1|
          \end{array}\right], \\
C_{+}^{(1)} & = &  \left[\begin{array}{cccc}
\lceil \frac{c_1}{R_1}\rceil R_1 - c_1 &                                         &        &            \\
                                   & (\lceil \frac{c_1}{R_1}\rceil + 1)R_1 - c_1 &        &            \\
                                   &                                             & \ddots &            \\
                                   &                                             &        & NR_1 - c_1
\end{array}\right], \\
C_{-}^{(2)} & = &  \left[\begin{array}{cccc}
|- c| &            &        & \\
          & |R_1 - c| &        & \\
          &            & \ddots & \\
          &            &        & |\lfloor \frac{c}{R_1} \rfloor R_1 - c|
          \end{array}\right], \\
C_{+}^{(2)} & = &  \left[\begin{array}{cccc}
\lceil \frac{c}{R_1}\rceil R_1 - c &                                       &        &            \\
                                   & (\lceil \frac{c}{R_1}\rceil + 1)R_1 - c &        &            \\
                                   &                                       & \ddots &            \\
                                   &                                       &        & NR_1 - c
\end{array}\right].
\end{eqnarray*}
We illustrate in Figure \ref{fig:infbuffer1}  the relationships between the large cast of characters.
\begin{figure}[b]
\begin{center}
\caption{Buffer~1} \label{fig:infbuffer1}
\begin{picture}(100,0)(0,0)
         \put(48,68){\makebox(0,0){{$T_1 = \left[\begin{array}{cc} T_{--}^{(1)} & T_{-+}^{(1)} \\
                                                             \vspace*{-0.2cm}  \\
                                                             T_{+-}^{(1)} & T_{++}^{(1)}
                                                             \end{array}\right]$}}}
         \put(33,43){\makebox(0,0){{$C_{-}^{(1)}, C_{+}^{(1)}$}}}
         \put(48,126){\makebox(0,0){{$T_1 = \left[\begin{array}{cc} T_{--}^{(2)} & T_{-+}^{(2)} \\
                                                              \vspace*{-0.2cm}  \\
                                                              T_{+-}^{(2)} & T_{++}^{(2)}
                                                              \end{array}\right]$}}}
         \put(33,100){\makebox(0,0){{$C_{-}^{(2)}, C_{+}^{(2)}$}}}
         \put(-15,90){\makebox(0,0){{$x^*$}}}
         \put(139,118){\makebox(0,0){{$\mathcal{S}_{-}^{(2)} \cup \mathcal{S}_{+}^{(2)}$}}}
         \put(152,90){\makebox(0,0){{$\mathcal{S}_{d}^{(*)} \cup \mathcal{S}_{s}^{(*)} \cup \mathcal{S}_{u}^{(*)}$}}}
         \put(-15,35){\makebox(0,0){{$0$}}}
         \put(139,62){\makebox(0,0){{$\mathcal{S}_{-}^{(1)} \cup \mathcal{S}_{+}^{(1)}$}}}
         \put(139,35){\makebox(0,0){{$\mathcal{S}_{s}^{(\circ)} \cup \mathcal{S}_{u}^{(\circ)}$}}}
\end{picture}
\end{center}
\end{figure}
Exploiting Markov-renewal arguments, da Silva Soares and Latouche \cite[Theorem~4.2]{leveldep} prove that the stationary density vector of a Markov-modulated level-dependent single-buffer fluid queue can be obtained by properly combining limiting densities from above and below each boundary (when possible) and steady state probability masses at these boundaries. To that effect, we consider the jump chain $\{J_n: n \geq 0 \}$ of the process $\{X(t), \varphi_1(t)\}$ restricted to the set of boundary states $\mathcal{B} = \{(\bullet, i): \bullet \in \{\circ, *\}, i \in \mathcal{S}\}$. We note that this jump chain will also be useful for obtaining bounds on marginal probabilities of Buffer~2, as described in Section~\ref{sec:infb2}. By \cite[Theorem~4.4]{leveldep}, the $(2N + 2) \times (2N + 2)$ transition matrix $\Omega$ of $\{J_n\}$, block-partitioned according to $\mathcal{B} = (\circ, \mathcal{S}_{u}^{(\circ)}) \cup (\circ, \mathcal{S}_s^{(\circ)}) \cup (*,\mathcal{S}_u^{(*)}) \cup (*,\mathcal{S}_s^{(*)}) \cup (*,\mathcal{S}_d^{(*)})$,  is
\begin{eqnarray}
\Omega & = & \left[\begin{array}{cc|ccc}
\cdot        & \Psi_{us}^{(\circ)}          & \Lambda_{uu}^{(\circ, *)}    & \Lambda_{us}^{(\circ,*)}    & \cdot \\
P_{su}^{(\circ)} & P_{ss}^{(\circ)}             & \cdot                 &              \cdot    & \cdot \\
\hline
&&&&\vspace*{-0.3cm} \\
\cdot        & \cdot                    & \cdot                 & \Psi_{us}^{(*)}       & \Psi_{ud}^{(*)}  \\
\cdot        & \cdot                    & P_{su}^{(*)}          & P_{ss}^{(*)}          & P_{sd}^{(*)}     \\
\cdot        & \hat{\Lambda}^{(*,\circ)}_{ds} & \hat{\Psi}^{(*)}_{du} & \hat{\Psi}^{(*)}_{ds} & \cdot
\end{array}\right], \label{eqn:Omega}
\end{eqnarray}
where $\Psi_{u m}^{(\bullet)}, \hat{\Psi}_{d m}^{(*)}, \Lambda_{u m}^{(\circ, *)}, \hat{\Lambda}_{d m}^{(*, \circ)}$ and $P_{s m}^{(\bullet)}$ denote various \emph{first passage probability} matrices, with \\
\hspace*{0.5cm} $[\Psi_{u m}^{(\bullet)}]_{ij}$ = probability of returning to $\bullet$ and in $j \in \mathcal{S}_{m}^{(\bullet)}$, after initially increasing from $\bullet$ and in $i \in \mathcal{S}_{u}^{(\bullet)}$, while avoiding a higher boundary (if there exists one), \\
\hspace*{0.5cm} $[\hat{\Psi}_{d m}^{(*)}]_{ij}$ = probability of returning to $x^*$ and in $j \in \mathcal{S}_{m}^{(*)}$, after initially decreasing from $x^*$ and in $i \in \mathcal{S}_{d}^{(*)}$, while avoiding level 0, \\
\hspace*{0.5cm} $[\Lambda_{u m}^{(\circ, *)}]_{ij}$ = probability of reaching $x^*$ and in $j \in \mathcal{S}_{m}^{(*)}$, while avoiding level 0 after initially increasing from there in $i \in \mathcal{S}_{u}^{\circ}$, \\
\hspace*{0.5cm} $[\hat{\Lambda}_{ds}^{(*, \circ)}]_{ij}$ = probability of reaching level 0 and in $j \in \mathcal{S}_{s}^{(\circ)}$, while avoiding $x^*$ after initially decreasing from there in $i \in \mathcal{S}_{d}^{*}$, and \\
\hspace*{0.5cm} $[P_{s m}^{(\bullet)}]_{ij}$ = probability of going from $i \in \mathcal{S}_{s}^{(\bullet)}$ to $j \in \mathcal{S}_{m}^{(\bullet)}$, in one transition. \\
\\
The jump chain of the Markov process $\{\varphi_1(\cdot)\}$ has the transition matrix
\begin{eqnarray} \label{eqn:P}
P = I - \Delta^{-1} T_1,
\end{eqnarray}
where $\Delta$ is the diagonal matrix with $[\Delta]_{i} = [T_1]_{ii}$--- for the remainder of the paper, we denote by $I$ the identity matrix of the appropriate size.
Clearly each of $P_{su}^{(\circ)}, P_{ss}^{(\circ)},$ $P_{su}^{(*)}, P_{ss}^{(*)}$ and $P_{sd}^{(*)}$ is a sub-matrix of $P$:
\begin{eqnarray*}
P = \left[\begin{array}{cc}
P_{ss}^{(\circ)} & P_{su}^{(\circ)} \\
P_{uu}^{(\circ)} & P_{us}^{(\circ)}
\end{array}\right] = \left[\begin{array}{ccc}
P_{dd}^{(*)} & P_{ds}^{(*)} & P_{du}^{(*)} \\
\vspace*{-0.3cm} \\
P_{sd}^{(*)} & P_{ss}^{(*)} & P_{su}^{(*)} \\
\vspace*{-0.3cm} \\
P_{ud}^{(*)} & P_{us}^{(*)} & P_{uu}^{(*)}
\end{array}
\right].
\end{eqnarray*}
The matrices $\Psi_{us}^{(\circ)}, [\Lambda_{u s}^{(\circ, *)},\Lambda_{u u}^{(\circ, *)}], \hat{\Lambda}_{ds}^{(*, \circ)}$ and $[\hat{\Psi}_{ds}^{(*)},\hat{\Psi}_{du}^{(*)}]$ are respectively equal to $\Psi_{+-}^{(1)}, \Lambda_{++}^{(1)}, \hat{\Lambda}_{--}^{(1)}$ and $\hat{\Psi}_{-+}^{(1)}$, the corresponding first passage probability matrices for the level-independent fluid queue $\{M_1(t),\rho_1(t): t \in \mathds{R}^{+}\}$ with finite size $x^*$, state space $\mathcal{S}_{-}^{(1)} \cup \mathcal{S}_{+}^{(1)}$, generator $T_1$ and rate matrices $C_{-}^{(1)}$ and $C_{+}^{(1)}$. By \cite[Theorem~5.2]{anafinite},
\begin{eqnarray} \label{eqn:bytheorem52}
& \left[\begin{array}{cc} \Lambda_{++}^{(1)} & \Psi_{+-}^{(1)} \\
                         \hat{\Psi}_{-+}^{(1)} & \hat{\Lambda}_{--}^{(1)} \end{array}\right] =
 \left[\begin{array}{cc} e^{\hat{U}_1x^*} & \Psi_1 \\
                         \hat{\Psi}_1 & e^{U_1x^*} \end{array}\right] \left[\begin{array}{cc}I & \Psi_1 e^{U_1x^*} \\
                       \hat{\Psi}_1e^{\hat{U}_1x^*} & I \end{array}\right]^{-1},
\end{eqnarray}
where $\Psi_1$ is the minimum nonnegative solution to the Riccati equation
\begin{eqnarray}
& (C_{+}^{(1)})^{-1}T_{+-}^{(1)} + (C_{+}^{(1)})^{-1}T_{++}^{(1)}\Psi_1 + \Psi_1 (C_{-}^{(1)})^{-1}T_{--}^{(1)} + \Psi_1 (C_{-}^{(1)})^{-1}T_{-+}^{(1)} \Psi_1 = 0, \label{eqn:Psi1}
\end{eqnarray}
$\hat{\Psi}_1$ is the minimum nonnegative solution to the Riccati equation
\begin{eqnarray}
 (C_{-}^{(1)})^{-1}T_{-+}^{(1)} + (C_{-}^{(1)})^{-1}T_{--}^{(1)}\hat{\Psi}_1 + \hat{\Psi}_1 (C_{+}^{(1)})^{-1}T_{++}^{(1)} + \hat{\Psi}_1 (C_{+}^{(1)})^{-1}T_{+-}^{(1)} \hat{\Psi}_1 = 0, \label{eqn:hatPsi1}
 \end{eqnarray}
 \begin{eqnarray}
      U_1 & = & (C_{-}^{(1)})^{-1}T_{--}^{(1)} + (C_{-}^{(1)})^{-1}T_{-+}^{(1)} \Psi_1, \label{eqn:U1}
\end{eqnarray}
and
\begin{eqnarray}
\hat{U}_1 & = & (C_{+}^{(1)})^{-1}T_{++}^{(1)} + (C_{+}^{(1)})^{-1}T_{+-}^{(1)} \hat{\Psi}_1. \label{eqn:hatU1}
\end{eqnarray}
Similarly, $[\Psi_{us}^{(*)},\Psi_{ud}^{(*)}] = \Psi_2$, which is the first passage probability matrix for the infinite level-independent fluid queue $\{M_2(t),\rho_2(t): t \in \mathds{R}^{+}\}$ with state space $\mathcal{S}_{-}^{(2)} \cup \mathcal{S}_{+}^{(2)}$, the generator $T_1$ and rate matrices $C_{-}^{(2)}$ and $C_{+}^{(2)}$. By \cite{rogers}, the matrix $\Psi_2$ is the minimum nonnegative solution to the Riccati equation
\begin{eqnarray}
(C_{+}^{(2)})^{-1}T_{+-}^{(2)} + (C_{+}^{(2)})^{-1}T_{++}^{(2)}\Psi_2 + \Psi_2 (C_{-}^{(2)})^{-1}T_{--}^{(2)} + \Psi_2 (C_{-}^{(2)})^{-1}T_{-+}^{(2)} \Psi_2 = 0. \label{eqn:Psi2}
\end{eqnarray}
Applying fast convergent algorithms described in \cite{ana02, bean05}, we can solve Riccati equations $(\ref{eqn:Psi1}), (\ref{eqn:hatPsi1})$ and $(\ref{eqn:Psi2})$ to obtain $\Psi_1, \hat{\Psi}_1$ and $\Psi_2$, and consequently $\Omega$. \\
\\
We denote by $\underline{m}= [\underline{p}^{(\circ)}_{s}, \underline{p}^{(*)}_{s}]$ the probability mass vector of Buffer~1 at the set of boundary sticky states $\mathcal{K} = \{(\bullet, \zeta): \bullet \in \{\circ, *\}, \zeta \in \mathcal{S}^{(\bullet)}_{s} \}$, and define $\mathcal{E}^{(*)} = \{(*,\zeta): \zeta \in \mathcal{S}_{u}^{(*)} \cup \mathcal{S}_d^{(*)}\} $ and $\mathcal{E}^{(\circ)} = \{(\circ,\zeta): \zeta \in \mathcal{S}_{u}^{(\circ)}\}$. Note that $\mathcal{K} = \mathcal{B} - \{\mathcal{E}^{(*)} \cup \mathcal{E}^{(\circ)}\}$. Proceeding in two steps, we write the transition matrix $\Omega^{(*)}$ of the censored fluid queue on $\{\mathcal{B} - \mathcal{E}^{(*)}\}$ as
\begin{eqnarray}
\Omega^{(*)} = \left[\begin{array}{ccc}
           \cdot & \Psi_{us}^{(\circ)} & \Lambda_{us}^{(\circ,*)} \\
P_{su}^{(\circ)} & P_{ss}^{(\circ)}    & \cdot \\
\cdot            & \cdot               & P_{ss}^{(*)}
\end{array}\right] + \left[\begin{array}{cc}
\Lambda_{uu}^{(\circ,*)} & \cdot \\
                   \cdot & \cdot \\
                   P_{su}^{(*)} & P_{sd}^{(*)}
\end{array}\right]\left[\begin{array}{cc}
                I & -\Psi_{ud}^{(*)} \\
-\hat{\Psi}_{ud}^{(*)} & I
\end{array}\right]^{-1}\left[\begin{array}{ccc}
\cdot & \cdot & \Psi_{us}^{(*)} \\
\cdot & \hat{\Lambda}_{ds}^{(*,\circ)} & \hat{\Psi}_{ds}^{(*)}
\end{array}\right], \nonumber
\end{eqnarray}
and find that the transition matrix of the censored fluid queue on $\mathcal{K}$ is
\begin{eqnarray}
 \Omega^{(\circ)} = \Omega^{(*)}_{\mathcal{K}\mathcal{K}} + \Omega^{(*)}_{\mathcal{K}\mathcal{E}^{(\circ)}}\left\{I - \Omega^{(*)}_{\mathcal{E}^{(\circ)}\mathcal{E}^{(\circ)}}\right\}^{-1}\Omega^{(*)}_{\mathcal{E}^{(\circ)}\mathcal{K}},
\end{eqnarray}
and its generator matrix is
\begin{eqnarray*}
\Theta = \Delta^{(s)} (I - \Omega^{(\circ)}),
\end{eqnarray*} where $\Delta^{(s)}$ is the diagonal matrix with $[\Delta^{(s)}]_{ii} = [T_1]_{ii}$ for $0 \leq i \leq \lfloor\frac{c}{R_1}\rfloor$. By \cite[Theorems~4.5 and 4.2]{leveldep},
\begin{eqnarray}
\underline{m} & = \kappa [\underline{x}_s^{(\circ)}, \underline{x}_s^{(*)}],
\end{eqnarray}
{and the density vector $\underline{\pi}(x)$ of Buffer~1 is}
\begin{eqnarray}
\underline{\pi}(x) & = & \kappa \underline{y}_1(x) \quad \mbox{ for } 0 < x < x^*, \\
                    & = & \kappa \underline{y}_2(x) \quad \mbox{ for } x > x^*,
\end{eqnarray}
where
\begin{eqnarray}
\kappa & = \left\{[\underline{x}_s^{(\circ)}, \underline{x}_s^{(*)}]\underline{1} + \int_0^{x^*} \underline{y}_{1} (x) \underline{1} dx + \int_{x^*}^{\infty} \underline{y}_{2} (x) \underline{1} dx\right\}^{-1}, \nonumber
\end{eqnarray}
the vector $[\underline{x}_s^{(\circ)}, \underline{x}_s^{(*)}]$ is a solution of $[\underline{x}_s^{(\circ)}, \underline{x}_s^{(*)}]\Theta = \underline{0}$,
\begin{eqnarray}
\underline{y}_2(x) & = & \{\underline{u}C_{+}^{(2)}N_{+}^{(2)}(x - x^*)\}(C^{(2)})^{-1},  \\
\underline{y}_1(x) & = & \{\underline{x}_s^{(\circ)}T_{-+}^{(1)}N_{+}^{(1)}(0,x) + \underline{d}C_{-}^{(1)}N_{-}^{(1)}(x^*, x)\}(C^{(1)})^{-1},
\end{eqnarray}

{the vectors $\underline{u}$ and $\underline{d}$ are the solution of}
\begin{eqnarray}
\underline{d} & = & \{\underline{x}_s^{(*)}T_{s-}^{(2)} + \underline{u}C_{+}^{(2)}\Psi_2\}(C_{-}^{(1)})^{-1},  \label{eqn:d} \\
\underline{u} & = & \{\underline{x}_s^{(*)}T_{s+}^{(2)} + \underline{x}_{s}^{(\circ)}T_{-+}^{(1)}\Lambda_{uu}^{(\circ,*)} + \underline{d}C_{-}^{(1)} \hat{\Psi}^{(1)}_{-+}\}(C_+^{(2)})^{-1}, \label{eqn:u}
\end{eqnarray}
and $N^{(2)}_{+}(w)$ is the matrix of expected number of visits to $w > 0$ in a phase of $\mathcal{S}_{+}^{(2)}$, while avoiding 0 after initially increasing from there, for the infinite fluid queue $\{M_2(t), \rho_2(t)\}$, \\
\hspace*{0.5cm} $N^{(1)}_{-}(x^*, w)$ is the matrix of expected number of visits to $w < x^*$ in a phase of $\mathcal{S}_{-}^{(1)}$, after initially decreasing from $x^*$ and while avoiding both $x^*$ and 0, for the finite fluid queue $\{M_1(t), \rho_1(t)\}$, and \\
\hspace*{0.5cm} $N^{(1)}_{+}(0, w)$ is the matrix of expected number of visits to $w < x^*$ in a phase of $\mathcal{S}_{+}^{(1)}$, after initially increasing from 0 and while avoiding both $0$ and $x^*$, for $\{M_1(t), \rho_1(t)\}$.\\
\\
By \cite[Theorems~2.1 and 2.2]{ram99},
\begin{eqnarray}
N^{(2)}_{+}(w) & = e^{K_2 w},
\end{eqnarray}
{with}
\begin{eqnarray}
 K_2 & = (C_{+}^{(2)})^{-1}T_{++}^{(2)} + \Psi_2 (C_{-}^{(2)})^{-1} T_{-+}^{(2)}, \nonumber
\end{eqnarray}
and by \cite[Lemma~4.1]{anafinite},
\begin{eqnarray}
\left[\begin{array}{l}
N_{+}^{(1)}(0,w) \\
\vspace*{-0.2cm} \\
N_{-}^{(1)}(x^*,w) \end{array}\right] = \left[\begin{array}{cc}
                                  I & e^{K_1x^*}\Psi_1 \\
      e^{\hat{K}_1 x^*}\hat{\Psi}_1 & I  \end{array}\right]^{-1}\left[\begin{array}{cc}
                                                                       e^{K_1w} & e^{K_1w}\Psi_1 \\
                                        e^{\hat{K}_1(x^* - w)}\hat{\Psi}_1  & e^{\hat{K}_1(x^* - w)}
                             \end{array}\right], \nonumber
\end{eqnarray}
with
\begin{eqnarray}
     K_1 & = (C_{+}^{(1)})^{-1}T_{++}^{(1)} + \Psi_1 (C_{-}^{(1)})^{-1} T_{-+}^{(1)}, \nonumber \\
\hat{K}_1 & = (C_{-}^{(1)})^{-1}T_{--}^{(1)} + \hat{\Psi}_1(C_{+}^{(1)})^{-1}T_{+-}^{(1)}. \nonumber
\end{eqnarray}
%

Therefore, the marginal distribution function of Buffer~1 is
\begin{eqnarray*}
\lim_{t \rightarrow \infty} P(X(t) \leq x) & = & \underline{p}_{s}^{(\circ)}\underline{1} + \int_{0}^{x} \underline{\pi}(x) dx \quad \mbox{ for } 0 < x < x^*, \\
& = & [\underline{p}_{s}^{(\circ)}, \underline{p}_{s}^{(*)}]\underline{1} + \int_{0}^{x} \underline{\pi}(x) dx \quad \mbox{ for } x \geq x^*.
\end{eqnarray*}
\section{Analysis for Buffer~2} \label{sec:infb2}

To derive the marginal probability distribution for Buffer~2 is not easy. Since its output capacity is dependent on $X(t)$, when analyzed as a standalone process, $\{Y(t), \varphi_2(t): t \in \mathds{R}^{+}\}$, Buffer~2 does not enjoy the Markovian property of $\{X(t), \varphi_1(t): t \in \mathds{R}^{+}\}$. Gautam \emph{et al.} \cite{gautam99} give bounds for the stationary distribution of fluid models with semi-Markov inputs and constant outputs. To apply these results, we first need to transform Buffer~2 into an equivalent fluid queue with semi-Markov inputs and a constant output. We achieve the transformation by employing a \emph{compensating source}, a concept developed by Elwalid and Mitra~\cite{elwalid95} and extended in Mahabhashyam \emph{et al.} \cite{resourcesharing}. The role of a compensating source is to add the exact amount of input for maintaining a constant output, $c$ in our case, while keeping all the time the fluid level the same as that of the original, output-varying, buffer. \\
\\
Consider a virtual fluid queue $\{Z(t), A(t), \varphi_2(t): t \in \mathds{R}^{+}\}$ which has $N$ exponential ON-OFF inputs and one independent compensating source. Here, $Z(t) \geq 0$ is the level, $A(t)$ is the semi-Markov process that drives the compensating source, and $\varphi_2(t)$ is the irreducible Markov chain controlling ON-OFF inputs, with state space $\mathcal{S}$ and generator $T_2$. The semi-Markov process $A(t)$ has state space $\mathcal{B}$, the set of boundary states for the jump chain $\{J_n\}$ defined in Section~\ref{sec:infb1} for the analysis of Buffer~1, as the output capacity of Buffer~2, and consequently the compensating source, changes each time $X(t)$ is in a boundary state. Specifically, the input rates $\dot{a}_{\bullet, i}$ of the compensating source are
\begin{eqnarray*}
\begin{array}{rll}
\dot{a}_{\bullet, i}  = & i R_1 & \mbox{for } (\bullet, i) \in (\circ, \mathcal{S}_{s}^{(\circ)}) \cup (*,\mathcal{S}_{s}^{(*)}), \\
                                        = & c_1   & \mbox{for } (\bullet, i) \in (\circ, \mathcal{S}_{u}^{(\circ)}) \cup (*,\mathcal{S}_{d}^{(*)}), \\
                                       = & c   & \mbox{for } (\bullet, i) \in (*, \mathcal{S}_{u}^{(*)}).
\end{array}
\end{eqnarray*}
Let $S_n$ be the time of the $n$th jump epoch in $A(t)$, $B_n$ the state of $A(t)$ immediately after the $n$th jump, and $\Omega(t)$ the kernel of $A(t)$, where
\begin{eqnarray*}
[\Omega(t)]_{ij} = P(S_1 \leq t, B_1 = j|B_0 = i).
\end{eqnarray*}
It is clear that $\Omega(\infty) = \Omega,$ the transition matrix of the jump chain $\{J_n\}$, given by (\ref{eqn:Omega}). We denote by $\widetilde{\Omega}(s)$ the matrix of Laplace-Stieltjes transforms of $S_1$, and in general, by $\widetilde{D}(s)$ and $\bar{D}(s)$ the respective LST counterparts of the sub-matrices $D$ and $\hat{D}$ of $\Omega$. The matrices $\widetilde{P}_{su}^{(\circ)}(s), \widetilde{P}_{ss}^{(\circ)}(s), \widetilde{P}_{su}^{(*)}(s),$ $\widetilde{P}_{ss}^{(*)}(s),$ and $\widetilde{P}_{sd}^{(*)}(s)$ are sub-matrices of $\widetilde{P}(s)$, where
\begin{eqnarray}
\widetilde{P}(s) = (sI -\Delta)^{-1}(T_1 - \Delta).
\end{eqnarray}
To obtain the remaining sub-matrices of $\widetilde{\Omega}(s)$, we follow an analogous analysis to that described in Section~\ref{sec:infb1}.
The matrices $\widetilde{\Psi}_{us}^{(\circ)}(s), [\widetilde{\Lambda}_{u s}^{(\circ, *)}(s),\widetilde{\Lambda}_{u u}^{(\circ, *)}(s)], \bar{\Lambda}_{ds}^{(*, \circ)}(s)$ and $[\bar{\Psi}_{ds}^{(*)}(s),\bar{\Psi}_{du}^{(*)}(s)]$ are equal to $\widetilde{\Psi}_{+-}^{(1)}(s), \widetilde{\Lambda}_{++}^{(1)}(s), \bar{\Lambda}_{--}^{(1)}(s)$ and $\bar{\Psi}_{-+}^{(1)}(s)$, the corresponding matrices of the LSTs of \emph{first passage times} for $\{M_1(t),\rho_1(t)\}$. By \cite[Theorem~3]{boundedmodel}, for $s$ such that Re$(s) > 0$,
\begin{eqnarray} \label{eqn:bytheorem3}
\left[\begin{array}{cc} \widetilde{\Lambda}_{++}^{(1)}(s) & \widetilde{\Psi}_{+-}^{(1)}(s) \\
                         \bar{\Psi}_{-+}^{(1)}(s) & \bar{\Lambda}_{--}^{(1)}(s) \end{array}\right]
 = \left[\begin{array}{cc} e^{\bar{U}_1(s)x^*} & \widetilde{\Psi}_1(s) \\
                         \bar{\Psi}_1(s) & e^{\widetilde{U}_1(s)x^*} \end{array}\right] \left[\begin{array}{cc}I & \widetilde{\Psi}_1(s) e^{\widetilde{U}_1(s)x^*} \\
                       \bar{\Psi}_1(s) e^{\bar{U}_1(s)x^*}& I \end{array}\right]^{-1}, \nonumber
\end{eqnarray}
where $\widetilde{\Psi}_1(s)$ is the minimum nonnegative solution to the Riccati equation
\begin{eqnarray}
(C_{+}^{(1)})^{-1}(T_{+-}^{(1)} - sI) + (C_{+}^{(1)})^{-1}(T_{++}^{(1)} - sI) \widetilde{\Psi}_1(s) \hspace*{4cm} \nonumber \\
+ \widetilde{\Psi}_1(s)(C_{-}^{(1)})^{-1}(T_{--}^{(1)} - sI) + \widetilde{\Psi}_1(s)(C_{-}^{(1)})^{-1}(T_{-+}^{(1)} - sI) \widetilde{\Psi}_1(s) = 0, \label{eqn:sPsi1}
\end{eqnarray}
$\bar{\Psi}_1(s)$ is the minimum nonnegative solution to the Riccati equation
\begin{eqnarray}
(C_{-}^{(1)})^{-1}(T_{-+}^{(1)} - sI) + (C_{-}^{(1)})^{-1}(T_{--}^{(1)} - sI)\bar{\Psi}_1(s) \hspace*{4cm} \nonumber \\
+ \bar{\Psi}_1(s) (C_{+}^{(1)})^{-1}(T_{++}^{(1)} - sI) + \bar{\Psi}_1(s) (C_{+}^{(1)})^{-1}(T_{+-}^{(1)} - sI) \bar{\Psi}_1(s) = 0,  \label{eqn:shatPsi1}
\end{eqnarray}
\begin{eqnarray}
      \widetilde{U}_1(s) = (C_{-}^{(1)})^{-1}(T_{--}^{(1)} - sI) + (C_{-}^{(1)})^{-1}(T_{-+}^{(1)} - sI) \widetilde{\Psi}_1(s), \nonumber
\end{eqnarray}
and
\begin{eqnarray}
\bar{U}_1(s) & = & (C_{+}^{(1)})^{-1}(T_{++}^{(1)} - sI) + (C_{+}^{(1)})^{-1}(T_{+-}^{(1)} - sI) \bar{\Psi}_1(s). \nonumber
\end{eqnarray}
Similarly, $[\widetilde{\Psi}_{us}^{(*)}(s),\widetilde{\Psi}_{ud}^{(*)}(s)] = \widetilde{\Psi}_2(s)$, which is the matrix of the LST of first passage times for $\{M_2(t),\rho_2(t)\}$. By \cite[Theorem~1]{unboundedfluids}, $\widetilde{\Psi}_2(s)$ is the minimum nonnegative solution to the Riccati equation
\begin{eqnarray}
(C_{+}^{(2)})^{-1}(T_{+-}^{(2)} - sI) + (C_{+}^{(2)})^{-1}(T_{++}^{(2)} - sI) \widetilde{\Psi}_2(s) \hspace*{4cm} \nonumber \\
+ \widetilde{\Psi}_2 (s) (C_{-}^{(2)})^{-1}(T_{--}^{(2)} - sI) + \widetilde{\Psi}_2(s) (C_{-}^{(2)})^{-1}(T_{-+}^{(2)} - sI)\widetilde{\Psi}_2(s) = 0.  \label{eqn:sPsi2}
\end{eqnarray}
Bean \emph{et al.} \cite{algosforlst} give efficient algorithms for solving $(\ref{eqn:sPsi1}), (\ref{eqn:shatPsi1})$ and $(\ref{eqn:sPsi2})$ to obtain $\widetilde{\Psi}_1(s), \bar{\Psi}_1(s)$ and $\widetilde{\Psi}_2(s)$, and consequently $\widetilde{\Omega}(s)$. \\
\\
Before we state the bounds for Buffer~2, we need to define \emph{effective bandwidths} and \emph{failure rate functions}.
For $v > 0$, the \emph{effective bandwidth} $\textsf{eb}(v)$ of an input that generates $F(t)$ amount of fluid in time $t$ is defined to be
\begin{eqnarray*}
\textsf{eb}(v) = \lim_{t \rightarrow \infty} \frac{1}{vt}\log E[e^{v F(t)}],
\end{eqnarray*}
(see, for example, \cite{elwalid93, kelly96}). By \cite{anick82,elwalid93}, the effective bandwidth $\textsf{eb}_{e}(v)$ of a single exponential ON-OFF source for fixed $v$ is
\begin{eqnarray}
\textsf{eb}_{e}(v) = \frac{R_2v - \alpha_2 - \beta_2 + \sqrt{(R_2v - \alpha_2 - \beta_2)^2 + 4\beta_2R_2 v}}{2v}.
\end{eqnarray}
To obtain the effective bandwidth $\textsf{eb}_{c}(v)$ for the compensating source, we begin by defining
$\Phi(v,u)$ to be the matrix with sub-matrices
\begin{eqnarray}
[\Phi(v,u)]_{(\bullet, i), (\bar{\bullet}, i')} = [\widetilde{\Omega}(v(u - \dot{a}_{\bullet, i}))]_{(\bullet, i), (\bar{\bullet}, i')}.
\end{eqnarray}
Denote by $\chi(D)$ the maximal real eigenvalue of a matrix $D$, then by \cite[Sections~4 and 5]{gautam99}, the effective bandwidth $\textsf{eb}_{c}(v)$ for fixed $v$ is the unique positive solution to the equation
\begin{eqnarray} \label{eqn:eb}
\chi(\Phi(v,\textsf{eb}_{c}(v))) = 1.
\end{eqnarray}
%
With these, we define $\eta$ to be the minimum positive solution to
\begin{eqnarray} \label{eqn:eta}
\textsf{eb}_c(\eta) + N\textsf{eb}_{e}(\eta) = c.
\end{eqnarray}
The existence of such $\eta$ is guaranteed by the facts \cite[Section 2.2.2]{gautamthesis} that $\textsf{eb}_{c}(v)$ and $\textsf{eb}_{e}(v)$ are both increasing functions with respect to $v$ and that for any given $v > 0$,
\begin{eqnarray*}
0 \leq \textsf{eb}_c(v) \leq c,
\end{eqnarray*}
and
\begin{eqnarray*}
\lim_{v \rightarrow 0}  \textsf{eb}_c(v) = 0 \quad \mbox{ and } \quad \lim_{v \rightarrow \infty}  \textsf{eb}_c(v) = c.
\end{eqnarray*}
For fixed $v$, we can solve (\ref{eqn:eb}) using fixed point iteration, as $\chi(\Phi(v,u))$ is a decreasing function with respect to $u$ \cite[Section 2.2.3]{gautamthesis}, and solve (\ref{eqn:eta}) using bisection. \\
\\
For $i \in \mathcal{B}$, we denote by $\tau_i$ the expected sojourn time of $A(t)$ in $i$
\begin{eqnarray}
\tau_i & = -\sum_{j \in \mathcal{B}}[\widetilde{\Omega}'(0)]_{ij},
\end{eqnarray}
by $\underline{p}$ the vector with elements
\begin{eqnarray}
p_i & = \frac{\omega_i \tau_i}{\sum\limits_{j \in \mathcal{B}} \omega_i \tau_j},
\end{eqnarray}
where $\underline{\omega}$ is the stationary vector associated with $\Omega$ ($\underline{\omega}\Omega = \underline{1}$, $\underline{\omega}\underline{1} = 1$),
and by $\underline{h}$ the left eigenvector of $\Phi(\eta,\textsf{eb}_{c}(\eta))$ corresponding to the eigenvalue one. Now, we are ready to define
$\Xi_{\max}(i,j)$ and $\Xi_{\min}(i,j)$  as follows
\begin{eqnarray}
&  \Xi_{\min}(i,j)  = \frac{h_i\tau_i}{p_i}\inf_x f_{ij}(x), \nonumber
\\
&   \Xi_{\max}(i,j) =  \frac{h_i\tau_i}{p_i}\sup_x f_{ij}(x), \nonumber
\end{eqnarray}
where
\begin{eqnarray} \label{eqn:f}
f_{ij}(x) = \frac{\displaystyle \int_x^{\infty}e^{\eta(\dot{a}_i - \textsf{eb}_c(\eta))y}d[\Omega(y)]_{ij}}
{e^{\eta(\dot{a}_i - \textsf{eb}_c(\eta))x}\left\{ [\Omega]_{ij} - [\Omega(x)]_{ij}\right\}} .
\end{eqnarray}
Applying \cite[Theorems~6 and 7]{gautam99} and then simplifying using \cite[Section~4.2.4]{gautamthesis}, we obtain the following result.
\begin{theorem} For $x > 0$,
\begin{eqnarray}
K_*e^{-\eta x} \leq \lim_{t \rightarrow \infty} P(Y(t) > x) \leq K^*e^{-\eta x},
\end{eqnarray}
where
\begin{eqnarray}
&   K_* =  \displaystyle\frac{\left[\displaystyle\frac{R_2}{\emph{\textsf{eb}}_{e}(\eta)\alpha_2}\right]^{N}H_c}{\max\limits_{s,(i,j)}D(s)\Xi_{\max}(i,j)}, \nonumber
\end{eqnarray}
{and}
\begin{eqnarray}
&  K^* = \frac{\left[\displaystyle\frac{R_2}{\emph{\textsf{eb}}_{e}(\eta)\alpha_2}\right]^{N} H_c}{\min\limits_{s,(i,j)}D(s)\Xi_{\min}(i,j)}, \nonumber
\end{eqnarray}
with $i, j \in \mathcal{B}$ and $1 \leq s \leq N$ such that $\dot{a}_{i} + s R_2 > c$ and $[\Omega]_{i,j} > 0$, where
\begin{eqnarray*}
H_c &= & \sum_{i \in \mathcal{B}}\left[\displaystyle\frac{h_i}{\eta(\dot{a}_i - \emph{\textsf{eb}}_c(\eta))}\right]\left[\sum_{j \in \mathcal{B}} [\Phi(\eta,\textsf{\emph{eb}}_{c}(\eta))]_{ij} -1\right],  \\
D(s) &= & \left[\displaystyle\frac{\alpha_2 + \beta_2}{\alpha_2 \beta_2}\right]^{s}\left[\displaystyle\frac{(\alpha_2 + \beta_2)(R_2 - \emph{\textsf{eb}}_{e}(\eta))}{\textsf{\emph{eb}}_e(\eta)\alpha_2^2}\right]^{(N - s)}.
\end{eqnarray*}
\end{theorem}
For $i, j \in \mathcal{B}$, the \emph{failure rate function} $\lambda_{ij}(x)$ of the compensating source is
\begin{eqnarray}
\lambda_{ij}(x) & = \frac{[\Omega'(x)]_{ij}}{[\Omega]_{ij} - [\Omega(x)]_{ij}}.
\end{eqnarray}
The function $[\Omega(x)]_{ij}$ is said to be increasing (IFR) if $\lambda_{ij}(x)$ is an increasing function of $x$, and decreasing (DFR) if $\lambda_{ij}(x)$ is a decreasing function of $x$. In the cases where $[\Omega(x)]_{ij}$ is either IFR or DFR, $\Xi_{\max}(i,j)$ and $\Xi_{\min}(i,j)$ are given in Table~\ref{tab:valuesofPsi}. For a sticky state $i \in (\circ,\mathcal{S}_{s}^{(\circ)}) \cup (*,\mathcal{S}_s^{(*)})$, $[\Omega(x)]_{ij}$ has a constant failure rate $\lambda_{ij}$, and $\Xi_{\min}(i,j) = \Xi_{\max}(i,j)$. \\
\\
When $[\Omega(x)]_{ij}$ is neither IFR nor DFR, $\Xi_{\max}(ij)$ and $\Xi_{\min}(ij)$ may be estimated by numerical computation. The LST of the numerator of $f_{ij}(\cdot)$ is $-\widetilde{\Omega}(s - \eta(\dot{a}_i - \textsf{eb}_c(\eta)))$; hence, both the numerator and the denominator of $f_{ij}(x)$ are obtainable by numerical inversion of $\widetilde{\Omega}(\cdot)$.


\begin{table*}[t]
\center
\caption{$\Xi_{\max}(i,j)$ and $\Xi_{\min}(i,j)$ in simple cases.}
\begin{tabular}{|c|c|c|}
\hline
\vspace*{-0.2cm} & & \\
& $\Xi_{\max}(i,j)$  &  $\Xi_{\min}(i,j)$
\\
\vspace*{-0.2cm} & & \\
\hline
\vspace*{-0.2cm} & & \\
$\begin{array}{l}
\mbox{IFR,  $\dot{a}_i > \textsf{eb}_c(\eta)$, or} \\
\mbox{DFR, $\dot{a}_i \leq \textsf{eb}_c(\eta)$} \end{array} $
& $\displaystyle\frac{[\Phi(\eta,\textsf{eb}_{c}(\eta))]_{ij}\tau_i h_i}{[\Omega]_{ij}p_i}$ & $\displaystyle\frac{\tau_i h_i \lambda_{ij}(\infty)}{p_i(\lambda_{ij}(\infty) - \eta(\dot{a}_i - \textsf{eb}_{c}(\eta)))}$  \\
\vspace*{-0.2cm} & & \\
\hline
\vspace*{-0.2cm} & & \\
$\begin{array}{l}
\mbox{IFR, $\dot{a}_i \leq \textsf{eb}_c(\eta)$ or} \\
\mbox{DFR, $\dot{a}_i > \textsf{eb}_c(\eta)$}
\end{array}$ & $\displaystyle\frac{\tau_i h_i \lambda_{ij}(\infty)}{p_i(\lambda_{ij}(\infty) - \eta(\dot{a}_i - \textsf{eb}_{c}(\eta)))}$ &  $\displaystyle\frac{[\Phi(\eta,\textsf{eb}_{c}(\eta))]_{ij}\tau_i h_i}{[\Omega]_{ij}p_i}$ \\
\hline
\end{tabular}  \label{tab:valuesofPsi}
\end{table*}

\section{Finite Buffer~1, with one input} \label{sec:finb1}

In this section, we determine the marginal probability distribution of Buffer~1 in the particular case $N = 1$, with an added
assumption that it has finite size $V > x^*$. Our aim is to illustrate the difference in distributions of the finite Buffer~1, as derived here, and of the infinite Buffer~1, as in \cite{resourcesharing}. As mentioned in the Introduction, while the analysis in this section can be extended in a straightforward manner to the general case $N \geq 1$, we specifically consider the case $N = 1$ to better illustrate the analytic approach. We only carry out the analysis for Buffer~1, as the expressions for Buffer~2 remain the same.\\
\\
The assumptions of the reference model, stated in Section~\ref{sec:refmodel}, become $R_j > c$, for $j = 1, 2$, and $\beta_1(R_1 - c) < \alpha_1 - c$.
The imposed finiteness leads to a third boundary $X(t) =~V$, in addition to the two boundaries $X(t) = 0$ and $X(t) = x^*$, and the second band becomes
$x^* < X(t) < V$. All state spaces are simplified significantly:
\begin{eqnarray*}
& \mathcal{S}_{s}^{(\circ)} = \mathcal{S}_{-}^{(1)} = \mathcal{S}_{d}^{(*)} = \mathcal{S}_{-}^{(2)} = \{0\}, \\
& \mathcal{S}_{u}^{(\circ)} = \mathcal{S}_{+}^{(1)} = \mathcal{S}_{u}^{(*)} = \mathcal{S}_{+}^{(2)} = \{1\}.
\end{eqnarray*}
While the set $\mathcal{S}_{s}^{(*)}$ of sticky states at $x^*$ is empty, there is a new sticky state at $V$, that is, $\mathcal{S}_{s}^{(V)} = \{1\}$ and $\mathcal{S}_{d}^{(V)} = \{0\}$. The generator matrix $T_1$ is
\begin{eqnarray*}
T_1 = \left[\begin{array}{cc} T_{--}^{(1)} & T_{-+}^{(1)} \\
                                                             \vspace*{-0.2cm}  \\
                                                             T_{+-}^{(1)} & T_{++}^{(1)}
                                                             \end{array}\right] = \left[\begin{array}{cc} T_{--}^{(2)} & T_{-+}^{(2)} \\
                                                             \vspace*{-0.2cm}  \\
                                                             T_{+-}^{(2)} & T_{++}^{(2)}
                                                             \end{array}\right] = \left[\begin{array}{rr} -\beta_1 & \beta_1 \\
\alpha_1 & -\alpha_1 \end{array}
\right],
\end{eqnarray*}
and the rate matrices are now scalars:
\begin{eqnarray*}
C_{+}^{(1)} &= R_1 - c_1, \qquad C_{-}^{(1)} = c_1, \\
C_{+}^{(2)} &= R_1 - c, \qquad \;\;C_{-}^{(2)} = c.
\end{eqnarray*}
By \cite[Theorem 4.4]{leveldep}, the jump chain $\{J_n: n \geq 0\}$ of the process $\{X(t),\varphi_1(t)\}$ restricted to the set of boundary states $\mathcal{B} = \{(\bullet, i): \bullet \in \{\circ,*,V\}, i = \{1,2\}\}$ has transition matrix
\begin{eqnarray}
\Omega = \left[\begin{array}{cc|cc|cc}
\cdot            & \Psi_{us}^{(\circ)}            & \Lambda_{uu}^{(\circ, *)}    & \cdot & \cdot & \cdot\\
1 & \cdot              & \cdot                        & \cdot & \cdot & \cdot \\
\hline
&&&&&\vspace*{-0.3cm} \\
\cdot            & \cdot                          & \cdot                        & \Psi_{ud}^{(*)} & \Lambda_{us}^{(*,V)} & \cdot \\
\cdot            & \hat{\Lambda}^{(*,\circ)}_{ds} & \hat{\Psi}^{(*)}_{du}        & \cdot           & \cdot                & \cdot \\
&&&&&\vspace*{-0.3cm} \\
\hline
&&&&&\vspace*{-0.3cm} \\
\cdot & \cdot & \cdot & \cdot & \cdot & 1 \\
\cdot & \cdot & \cdot & \hat{\Lambda}_{dd}^{(V,*)} & \hat{\Psi}_{ds}^{(V)} & \cdot
\end{array}\right], \label{eqn:Omegafinite}
\end{eqnarray}
where $\Psi_{us}^{(\circ)}, \Lambda_{u u}^{(\circ, *)}, \hat{\Lambda}_{ds}^{(*, \circ)}$ and $\hat{\Psi}_{du}^{(*)}$ are the solutions of (\ref{eqn:bytheorem52}). Here, (\ref{eqn:Psi1}) and (\ref{eqn:hatPsi1}) reduce to scalar quadratic equations, from which one easily obtains the minimal solutions
\begin{eqnarray}
\Psi_1 = 1, \qquad \hat{\Psi}_1 = \frac{\beta_1(R_1 - c_1)}{\alpha_1 c_1}. \label{eqn:Psis}
\end{eqnarray}
Substituting (\ref{eqn:Psis}) into (\ref{eqn:U1}) and (\ref{eqn:hatU1}) leads to
\begin{eqnarray}
U_1 = 0, \qquad \hat{U}_1 = \frac{-\alpha_1 c_1 + \beta_1(R_1 - c_1)}{c_1(R_1 - c_1)}  \label{eqn:Us1}
\end{eqnarray}
Then, substituting (\ref{eqn:Us1}) into $(\ref{eqn:bytheorem52})$ gives us
\begin{eqnarray} \label{eqn:bytheorem52solved}
   \Lambda_{uu}^{(\circ,*)} & = & \Lambda_{++}^{(1)} = \frac{1 - \hat{\Psi}_1}{e^{-\hat{U}_1x^*} - \hat{\Psi}_1}, \\
   \Psi_{us}^{(\circ)} & = & \Psi_{+-}^{(1)} = \frac{1 - e^{\hat{U}_1x^*}}{1 - \hat{\Psi}_1e^{\hat{U}_1x^*}}, \\
   \hat{\Psi}_{du}^{(*)} & = & \hat{\Psi}_{-+}^{(1)} = \frac{\hat{\Psi}_1 - \hat{\Psi}_1e^{\hat{U}_1x^*}}{1 - \hat{\Psi}_1e^{\hat{U}_1x^*}}, \\
   \hat{\Lambda}_{ds}^{(*,\circ)} & = & \hat{\Lambda}_{--}^{(1)} = \frac{1 - \hat{\Psi}_1}{1 - \hat{\Psi}_1e^{\hat{U}_1 x^*}}. \label{eqn:bytheorem52solved4}
\end{eqnarray}
The second band now being finite, we follow the same steps as for the first band and find that the matrices $\Psi_{ud}^{(*)}$, $\Lambda_{us}^{(*,V)}$, $\hat{\Lambda}_{dd}^{(V,*)}$ and $\hat{\Psi}_{ds}^{(V)}$ are equal to $\Psi_{+-}^{(2)}, \Lambda_{++}^{(2)}, \hat{\Lambda}_{--}^{(2)}$ and $\hat{\Psi}_{-+}^{(2)}$, the corresponding first passage probability matrices for the level-independent fluid queue $\{M_2(t),\rho_2(t): t \in \mathds{R}^{+}\}$ with finite size $V - x^*$, state space $\mathcal{S}_{-}^{(2)} \cup \mathcal{S}_{+}^{(2)}$, generator $T_1$ and rates $C_{-}^{(2)}$ and $C_{+}^{(2)}$.
By \cite[Theorem~5.2]{anafinite},
\begin{eqnarray}\label{eqn:bytheorem52again}
\left[\begin{array}{cc} \Lambda_{++}^{(2)} & \Psi_{+-}^{(2)} \\
                         \hat{\Psi}_{-+}^{(2)} & \hat{\Lambda}_{--}^{(2)} \end{array}\right]
 =  \left[\begin{array}{cc} e^{\hat{U}_2(V - x^*)} & \Psi_2 \\
                         \hat{\Psi}_2 & e^{U_2(V - x^*)} \end{array}\right] \left[\begin{array}{cc}1 & \Psi_2 e^{U_2(V - x^*)} \\
                       \hat{\Psi}_2e^{\hat{U}_2(V - x^*)} & 1 \end{array}\right]^{-1}, \nonumber
\end{eqnarray}
where
\begin{eqnarray}
\Psi_2 & = & 1, \qquad \hat{\Psi}_2 = \frac{\beta_1(R_1 - c)}{\alpha_1 c}, \label{eqn:Psis2} \\
U_2 & = & 0, \qquad \hat{U}_2 = \frac{-\alpha_1 c + \beta_1(R_1 - c)}{c(R_1 - c)}. \label{eqn:Us2}
\end{eqnarray}
Substituting (\ref{eqn:Psis2}) and (\ref{eqn:Us2}) into $(\ref{eqn:bytheorem52again})$ gives us
\begin{eqnarray} \label{eqn:bytheorem52againsolved}
   \Lambda_{us}^{(*,V)} & = & \Lambda_{++}^{(2)} = \frac{1 - \hat{\Psi}_2}{e^{-\hat{U}_2(V - x^*)} - \hat{\Psi}_2}, \\
   \Psi_{ud}^{(*)} & = & \Psi_{+-}^{(2)} = \frac{1 - e^{\hat{U}_2(V - x^*)}}{1 - \hat{\Psi}_2e^{\hat{U}_2(V - x^*)}}, \\
   \hat{\Psi}_{ds}^{(V)} & = & \hat{\Psi}_{-+}^{(2)} = \frac{\hat{\Psi}_2 - \hat{\Psi}_2e^{\hat{U}_2(V - x^*)}}{1 - \hat{\Psi}_2e^{\hat{U}_2(V - x^*)}}, \\
   \hat{\Lambda}_{dd}^{(V,*)} & = & \hat{\Lambda}_{--}^{(2)} = \frac{1 - \hat{\Psi}_2}{1 - \hat{\Psi}_2e^{\hat{U}_2(V - x^*)}}. \label{eqn:bytheorem52againsolved4}
\end{eqnarray}
Together, equations (\ref{eqn:bytheorem52solved})--(\ref{eqn:bytheorem52solved4}) and (\ref{eqn:bytheorem52againsolved})--(\ref{eqn:bytheorem52againsolved4}) complete the transition matrix $\Omega$, specified in (\ref{eqn:Omegafinite}), of the jump chain $\{J_n\}$ on the set $\mathcal{B}$ of boundary states. The set $\mathcal{K}$ of sticky states is $\{(\circ,0), (V,1)\}$.
Straightforward but tedious calculations show that the jump chain on $\mathcal{K}$ has the transition matrix
\begin{eqnarray*}
\Omega^{(\circ)} = \left[\begin{array}{cc}
\Psi_{us}^{(\circ)} + \displaystyle\frac{\Lambda_{uu}^{(\circ,*)}\Psi_{ud}^{(*)}\hat{\Lambda}^{(*,\circ)}_{ds}}{1 - \Psi_{ud}^{(*)}\hat{\Psi}_{du}^{(*)}} & \displaystyle\frac{\Lambda_{uu}^{(\circ,*)}\Lambda_{us}^{(*,\circ)} }{1 - \Psi_{ud}^{(*)}\hat{\Psi}_{du}^{(*)}} \\
\vspace*{-0.2cm} \\
\displaystyle\frac{\hat{\Lambda}_{dd}^{(V,*)}\hat{\Lambda}_{ds}^{(*,\circ)}}{1 - \Psi_{ud}^{(*)}\hat{\Psi}_{du}^{(*)}} &  \hat{\Psi}_{ds}^{(V)} + \displaystyle\frac{\hat{\Lambda}_{dd}^{(V,*)}\hat{\Psi}_{du}^{(*)}\Lambda_{us}^{(*,V)}}{1 - \Psi_{ud}^{(*)}\hat{\Psi}_{du}^{(*)}}
\end{array}\right], \end{eqnarray*}
and, consequently, the generator matrix
\begin{eqnarray*}
 \Theta =  \left[\begin{array}{cc} -\beta_1 & \\
         & -\alpha_1 \end{array} \right](I - \Omega^{(\circ)}). \end{eqnarray*}
A solution of $[{x}_s^{(\circ)},{x}_s^{(V)}]\Theta = \underline{0}$ is
\begin{eqnarray*}
{x}_s^{(\circ)} & = & 1, \\
{x}_s^{(V)}      & = & -\frac{\beta_1(1 - [\Omega^{(\circ)}]_{11})}{\alpha_1[\Omega^{(\circ)}]_{21}}.
\end{eqnarray*}
By \cite[Theorem 4.5]{leveldep}, the probability mass vector $\underline{m} = [p_{s}^{(\circ)},{p}_{s}^{(V)}]$ of Buffer~1 at $\mathcal{K}$ is given by $\underline{m} = \kappa[x_{s}^{(\circ)},x_{s}^{(*)}]$, with
\begin{eqnarray*}
\kappa & = & \{1 + \frac{\alpha_1(1 - [\Omega^{(\circ)}]_{11})}{\beta_1[\Omega^{(\circ)}]_{22}}  + \int_0^{x^*} \underline{y}_{1} (x) \underline{1} dx  + \int_{x^*}^{V} \underline{y}_{2} (x) \underline{1} dx \}^{-1},
\end{eqnarray*}
and
\begin{eqnarray*}
\underline{y}_1(x) 
& = &  \{\beta_1N_{+}^{(1)}(0,x) + c_1\gamma_1N_{-}^{(1)}(x^*,x)\}(C^{(1)})^{-1}, \\
\underline{y}_2(x) 
& = & \{(R_1 - c)\gamma_2N_{+}^{(2)}(0, x - x^*) + c\gamma_3N_{-}^{(2)}(V - x^*, x - x^*)\}(C^{(2)})^{-1};
\end{eqnarray*}
the vectors $\gamma_1$, $\gamma_2$ and $\gamma_3$ are the solution of the system
\begin{eqnarray*}
\gamma_3 & = & {x}_s^{(V)}T_{+-}^{(2)}(C_{-}^{(2)})^{-1} \\
& = & \frac{\alpha_1}{c}{x}_s^{(V)}, \\
\gamma_2  & = & \{{x}_{s}^{(\circ)}T_{-+}^{(1)}\Lambda_{++}^{(\circ,*)} + \gamma_1C_{-}^{(1)}\hat{\Psi}_{-+}^{(1)}\}(C_{+}^{(2)})^{-1} \\
& = & \frac{1}{R_1 - c}\{\beta_1\Lambda_{++}^{(\circ,*)} + c_1\gamma_1\hat{\Psi}_{-+}^{(1)}\}, \\
\gamma_1 & = & \{\gamma_2C_{+}^{(2)}\Psi_{+-}^{(2)} + \gamma_3C_{-}^{(2)}\hat{\Lambda}_{--}^{(V,*)}\}(C_{-}^{(1)})^{-1} \\
& = & \frac{1}{c_1} \{(R_1 - c)\gamma_2\Psi_{+-}^{(2)} + \alpha_1{x}^{(V)}_{s}\hat{\Lambda}_{--}^{(V,*)}\}.
\end{eqnarray*}
Solving for $\gamma_1$ and $\gamma_2$ leads to
\begin{eqnarray*}
 \left[\begin{array}{c}
\gamma_1 \\
\gamma_2
\end{array}\right]  & =  & \left[\begin{array}{cc}
\displaystyle\frac{c_1\hat{\Psi}_{-+}^{(1)}}{R_1 - c}  & -1 \\
1 & -\displaystyle\frac{(R_1 - c)\Psi_{+-}^{(2)}}{c_1}
\end{array}\right]^{-1}\left[\begin{array}{c}
-\displaystyle\frac{\beta_1\Lambda_{++}^{(\circ,*)}}{R_1 - c}  \\
\vspace{-0.2cm} \\
\displaystyle\frac{\alpha_1{x}_s^{(V)}\hat{\Lambda}_{--}^{(V,*)}}{c_1}
\end{array}
\right] \\
& =  & \frac{1}{1-\hat{\Psi}_{-+}^{(1)}\Psi_{+-}^{(2)}}  \left[
\begin{array}{c}
\displaystyle\frac{\beta_1}{c_1}\Psi_{+-}^{(2)}\Lambda_{++}^{(\circ,*)} + \displaystyle\frac{\alpha_1{x}_s^{(V)}\hat{\Lambda}_{--}^{(V,*)}}{c_1}  \\
\vspace*{-0.2cm} \\
\displaystyle\frac{\beta_1}{R_1 - c}\Lambda_{++}^{(\circ,*)} +  \displaystyle\frac{\alpha_1{x}_s^{(V)}}{R_1 - c}\hat{\Psi}_{-+}^{(1)}\hat{\Lambda}_{--}^{(V,*)}
\end{array}
\right].
\end{eqnarray*}
By \cite[Lemma 4.1]{anafinite},
\begin{eqnarray}
\left[\begin{array}{l}
N_{+}^{(1)}(0,x) \\
\vspace*{-0.2cm} \\
N_{-}^{(1)}(x^*,x) \end{array}\right] & =  \left[\begin{array}{cc}
                                  1 & e^{K_1x^*} \\
      e^{\hat{K}_1 x^*}\hat{\Psi}_1 & 1  \end{array}\right]^{-1}\left[\begin{array}{cc}
                                                                       e^{K_1x} & e^{K_1x} \\
                                        e^{\hat{K}_1(x^* - x)}\hat{\Psi}_1  & e^{\hat{K}_1(x^* - x)}
                             \end{array}\right],
\end{eqnarray}
with
\begin{eqnarray}
     K_1 
             & = -\frac{\alpha_1}{R_1 - c_1} + \frac{\beta_1}{c_1}, \qquad 
\hat{K}_1 
                = 0. \nonumber
\end{eqnarray}
Consequently,
\begin{eqnarray*}
  N_{+}^{(1)}(0,x) & = \frac{1}{1 - \hat{\Psi}_1e^{K_1x^*}}\left[e^{K_1x} - \hat{\Psi}_1e^{K_1x^*}, e^{K_1x} - e^{K_1x^*}\right],
\end{eqnarray*}
 {and}
\begin{eqnarray*}
N_{-}^{(1)}(x^*,x) & = \frac{1}{1 - \hat{\Psi}_1e^{K_1x^*}}\left[-\hat{\Psi}_{1}e^{K_1x} + \hat{\Psi}_1, -\hat{\Psi}_{1}e^{K_1x} + 1\right].
\end{eqnarray*}
Similarly, by \cite[Lemma 4.1]{anafinite} again,
\begin{eqnarray}
 \left[\begin{array}{l}
N_{+}^{(2)}(0,x - x^*) \\
\vspace*{-0.2cm} \\
N_{-}^{(2)}(V - x^*,x - x^*) \end{array}\right]  =   \left[\begin{array}{cc}
                                  1 & e^{K_2(V - x^*)} \\
      e^{\hat{K}_2(V - x^*)}\hat{\Psi}_2 & 1  \end{array}\right]^{-1}\left[\begin{array}{cc}
                                                                       e^{K_2(x - x^*)} & e^{K_2(x - x^*)} \\
                                        e^{\hat{K}_2(V - x)}\hat{\Psi}_2  & e^{\hat{K}_2(V - x)}
                             \end{array}\right], \nonumber
   \end{eqnarray}
with
\begin{eqnarray}
     K_2 
             & = -\frac{\alpha_1}{R_1 - c} + \frac{\beta_1}{c}, \qquad 
\hat{K}_2 
                = 0. \nonumber
\end{eqnarray}
Consequently,
\begin{eqnarray*}
 N_{+}^{(2)}(0,x - x^*)
  = \frac{e^{-K_2x^*}}{1 - \hat{\Psi}_2e^{K_2(V - x^*))}}\left[e^{K_2x} - \hat{\Psi}_2e^{K_2V}, e^{K_2x} - e^{K_2V}\right],
 \end{eqnarray*}
and
\begin{eqnarray*}
 N_{-}^{(2)}(V - x^*,x - x^*)  = \frac{1}{1 - \hat{\Psi}_2e^{K_2(V - x^*)}}\left[-\hat{\Psi}_{2}e^{K_2(x - x^*)} + \hat{\Psi}_2, -\hat{\Psi}_{2}e^{K_2(x - x^*)} + 1\right].
\end{eqnarray*}
The density vector $\underline{\pi}(x)$ of Buffer~1 is
\begin{eqnarray*}
\underline{\pi}(x) & = & \kappa \underline{y}_1(x) \quad \mbox{ for } 0 < x < x^*, \\
                              & = & \kappa \underline{y}_2(x) \quad \mbox{ for } x^*< x < V.
\end{eqnarray*}
As an illustration, we consider Scenarios A, E, and F from \cite[Table~1]{resourcesharing}, to compare marginal probabilities for Buffer~1 in the finite and infinite cases. In all three scenarios,   $R_1 = 12.48,  \alpha_1 = 11, \beta_1 = 1, x^* = 1.5$, and $c = 2.6$. For Scenario A, $c_1 = 1.6$ and $c_2 = 1$; for Scenario B, $c_1 = 1.19$ and $c_2 = 1.41$; and, for Scenario C, $c_1 = 0.2$ and $c_2 = 2.4$.  In Tables \ref{tab:Pxstar} and \ref{tab:P3}, the values of $\lim_{t \rightarrow \infty} P(X(t) > x^*)$ and of $\lim_{t \rightarrow \infty}  P(X(t) > 3)$ for $V = \infty$ are taken from the last and the second columns of \cite[Table~2]{resourcesharing}, respectively. \\
\begin{table}[h]
\begin{center}
\caption{$\lim_{t \rightarrow \infty}  P(X(t) > x^*)$.} \label{tab:Pxstar}
\begin{tabular}{|c|l|l|l|l|}
\hline
    $V$             & $\infty$ &        $3.5$   & $ 6$              & $20$            \\
\hline
              A      & 0.1706   &  0.1411       &    0.1660       &    0.1706     \\
\hline
              E      & 0.1942   &  0.1615       &    0.1891       &   0.1942      \\
\hline
              F      & 0.3501   &  0.3009       &    0.3426       &   0.3501      \\
\hline
\end{tabular}
\end{center}
\end{table} \\
\begin{table}[h]
\begin{center}
\caption{$\lim_{t \rightarrow \infty}  P(X(t) > 3)$.} \label{tab:P3}
\begin{tabular}{|c|l|l|l|l|}
\hline
    $V$             & $\infty$  &        $3.5$   & $ 6$          & $20$            \\
\hline
              A       & 0.0572   &       0.0237 &   0.0519    & 0.0572         \\
\hline
              E       & 0.0651   &       0.0271  &  0.0592     & 0.0651        \\
\hline
              F       & 0.1173   &      0.0505   &  0.1072     & 0.1173        \\
\hline
\end{tabular}
\end{center}
\end{table}\\
It is clear that the marginal probabilities for Buffer~1 differ between the infinite and finite cases, and that these differences quickly tend to zero as
$V$ tends to infinity.


\section{Acknowledgement}
This work has been subsidized by the ARC grant AUWB-08/13-ULB 5 financed by the Minist\`{e}re de la Communaut\'{e} fran\c{c}aise de Belgique. The second author also gratefully acknowledges the hospitality of the Mathematical Institute at University of Wroc{\l}aw, where part of this work has been done.

\bibliographystyle{abbrv}
\bibliography{List}
\end{document}